\documentclass[12pt,cd]{amsart}
\usepackage{amsfonts, amssymb,amsmath, amsthm, amsxtra, latexsym, mathrsfs, amscd}
\usepackage{amsmath, amssymb, amsthm, times, float}
\usepackage{graphics}
\usepackage{pb-diagram}
\usepackage[latin1]{inputenc} %

\newtheorem{theo+}{Theorem}[section]
\newtheorem{prop+}[theo+]{Proposition}
\newtheorem{coro+}[theo+]{Corollary}
\newtheorem{lemm+} [theo+]{Lemma}
\newtheorem{deep+}  [theo+]  {Deep Result}
\newtheorem{fact+}  [theo+]  {Fact}
\theoremstyle{definition}
\newtheorem{exam+}  [theo+]  {Example}
\newtheorem{rema+}  [theo+]  {Remark}
\newtheorem{defi+}  [theo+]  {Definition}
\newtheorem{xca+}[theo+]{Exercise}

\numberwithin{equation}{section}

\usepackage{color}

\def\draft{\centerline{(Draft {\the \day}/{\the\month} \the \year.)}}
\bigskip

\def\refn#1.#2{\expandafter\def\csname#1\endcsname{[#2]}}
\def\refnr#1.{\csname#1\endcsname}

\def\fa{\mathfrak a}

\def\fg{\mathfrak g}
\def\fk{\mathfrak k}

\def\fl{\mathfrak l}
\def\fm{\mathfrak m}
\def\fn{\mathfrak n}
\def\fp{\mathfrak p}

\def\a{\alpha}

\def\Claminv2{|C(\Lambda)|^{-2}}

\def\Ga{\Gamma}
\def\varepsi{\varepsilon}

\def\de{d\varepsilon}

\def\Aa2D{A^{\a,2}(D)}
\def\bAa2D{\overline{A^{\a,2}(D)}}
\def\Ab2D{A^{\beta,2}(D)}
\def\bAb2D{\overline{A^{\beta,2}(D)}}

\def\Norm#1_#2{\Vert#1\Vert_{#2}}

\def\phipl12{\phi_{p_{l_1}, p_{l_2}}}
\def\phip01{\phi_{p_{0}, p_{0}}}

\def\a{\alpha}

\def\Claminv2{|C(\Lambda)|^{-2}}

\def\varepsi{\varepsilon}

\def\Ga{\Gamma}

\def\Ad{\operatorname{Ad}}

\def\det{\operatorname{det}}

\def\SL{\operatorname{SL}}

\def\de{d\varepsilon}

\def\Aa2D{A^{\a,2}(D)}
\def\bAa2D{\overline{A^{\a,2}(D)}}
\def\Ab2D{A^{\beta,2}(D)}
\def\bAb2D{\overline{A^{\beta,2}(D)}}

\def\phipl12{\phi_{p_{l_1}, p_{l_2}}}
\def\phip01{\phi_{p_{0}, p_{0}}}

\def\bc{\mathbb C}
\def\br{\mathbb R}

\def\bh{\mathbb H}

\def\alg/{algebra}

\def\Alg/{Algebra} 

\def\alt/{alternative} 
\def\anal/{analytic}
\def\analfunc/{\anal/\ \func/}
\def\Ans/{\it Answer. \normal}
\def\ass/{associative}
\def\nass/{non-\ass/}
\def\autom/{automorphism}
\def\homom/{homomorphism}
\def\isom/{isomorphism}
\def\bdd/{bounded}
\def\Bdd/{Bounded}
\def\bddsymdom/{bounded \sym/ \dom/}
\def\Cartdom/{Cartan \dom/}
\def\bdry/{boundary}
\def\bsd/{\bdd/ \symdom/}
\def\bv/{boundary value}
\def\cf/{{\it cf}\.}
\def\Cf/{{\it Cf}\.}
\def\charr/{character}
\def\coeff/{coefficient}
\def\comm/{commutative}
\def\cpct/{compact}
\def\compl/{complex}
\def\comp/{complex}
\def\Comp/{Complex}
\def\conf/{conformal}
\def\conj/{conjugate}
\def\conn/{connect}
\def\cont/{continuous}
\def\conv/{converge} 
\def\convc/{convergence}
\def\convt/{convergent}
\def\convx/{convex}
\def\coord/{coordinate}
\def\lcoord/{local coordinate}
\def\Corr/{Corresponding}
\def\corr/{corresponding}
\def\corrd/{correspond}
\def\cov/{covariant}
\def\decomp/{decomposition}
\def\deco/{decompose}
\def\diff/{different} 
\def\Diff/{Different} 
\def\dimn/{dimension} 
\def\distr/{distribution} 
\def\div/{diverge} 
\def\dom/{domain}
\def\eg/{\hbox{\it e.g}\.}
\def\eigenf/{eigen\-\func/}
\def\eigensp/{eigen\-space}
\def\eigenv/{eigen\-value}
\def\eq/{equation}
\def\equa/{equation}
\def\de/{\diff/ial \equa/}
\def\do/{\diff/ial operator}
\def\ode/{ordinary \de/}
\def\pde/{partial \de/}
\def\pdo/{partial \diff/ial operator}
\def\psdo/{pseudo \diff/ial operator}
\def\fin/{finite}
\def\Ex/{\it Example.\ \normal}
\def\Exnr#1/{\it Example #1.\ \normal}
\def\foll/{follow}
\def\follg/{following}
\def\Follg/{Following}
\def\func/{function}
\def\Func/{Function}
\def\Fonc/{Fonc\-tion}
\def\fonc/{fonc\-tion}
\def\Funk/{Funk\-tion}
\def\funk/{Funk\-tion}
\def\gen/{general}
\def\har/{harmonic}
\def\Hint/{\it Hint. \normal}
\def\hist/{historic}
\def\histcl/{historical}
\def\hol/{holo\-morphic}
\def\homog/{ho\-mo\-ge\-ne\-ous}
\def\hyp/{hyper\-bolic}
\def\hyperg/{hyper\-geometric}
\def\ie/{\hbox{\it i.e.}}
\def\iff/{if and only if}
\def\ineq/{inequality}
\def\infra/{{\it inf\-ra}}
\def\ultra/{{\it ult\-ra}}
\def\Inpart/{In particular}
\def\inpart/{in particular}
\def\instof/{instead of}
\def\interps/{interpolation space}
\def\interp/{interpolation}
\def\Interp/{Interpolation}
\def\interpr/{Interpretation}
\def\Intr/{Introduction}
\def\intv/{interval}
\def\inv/{invariant}
\def\invc/{invariance}
\def\Iowords/{In other words}
\def\iowords/{in other words}
\def\ipr/{inner product}
\def\irred/{irreducible}
\def\lb/{line bundle}
\def\lin/{linear}
\def\lhs/{left hand side}
\def\rhs/{right hand side}
\def\loc/{local}
\def\math/{mathematic} 
\def\mathcn/{\math/ian}
\def\manif/{manifold}
\def\meas/{measure}
\def\measl/{measurable}
\def\mero/{mero\-morphic}
\def\mon/{monomial}
\def\monog/{monogenic}
\def\mult/{multiple}
\def\multy/{multiply}
\def\multn/{multiplication}
\def\nas/{necessary and sufficient}
\def\nbd/{neighborhood}
\def\neg/{negative}
\def\nondeg/{nondegenerate}
\def\Oohand/{On the other hand}
\def\oohand/{on the other hand}
\def\Oonhand/{On the one hand}
\def\oonhand/{on the one hand}
\def\oper/{operator}
\def\orth/{ortho\-gonal}
\def\orthon/{ortho\-normal}
\def\otoh/{on the other hand}
\def\quat/{quaternion}
\def\pp/{\hbox{a. e.}}
\def\psh/{plurisubharmonic}
\def\pol/{polynomial}
\def\pot/{potential}
\def\pos/{positive}
\def\princ/{principle}
\def\prob/{probability}
\def\proj/{projective}
\def\projn/{projection}
\def\Proof/{\it Proof:\normal}
\def\Rem/{\it Remark\normal}
\def\Remnr#1/{\it Remark\ \normal #1. }
\def\rep/{representation}
\def\reps/{representations}
\def\meta/{metaplectic representation}
\def\repr/{reproducing}
\def\reprker/{reproducing kernel}
\def\resp/{respective} 
\def\resply/{respectively}
\def\restr/{restriction}
\def\sa/{self-adjoint}
\def\st/{such that}

\def\sol/{solution}
\def\ru/{space}
\def\sph/{spherical}
\def\ssp/{sub\ru/}
\def\sym/{symmetric}
\def\Sym/{Symmetric}
\def\symb/{symbol}
\def\symbc/{symbolic}
\def\symdom/{\sym/ domain}
\def\symp/{symplectic}
\def\Theor#1/{\fet Theorem #1.\ \normal}
\def\Lem#1/{\fet Lemma #1.\ \normal}
\def\Lemma/{\fet Lemma.\ \normal}
\def\topl/{topology}
\def\topll/{topological}
\def\transf/{transform}
\def\transl/{translation}
\def\transfn/{transformation}
\def\transv/{transvectant}
\def\trig/{trigonometric}
\def\tril/{trilinear}
\def\trilf/{trilinear form}
\def\uhp/{upper halfplane}
\def\uhs/{upper halfspace}
\def\vb/{vector bundle}
\def\vf/{vector field}
\def\vsp/{vector space}
\def\wrt/{with respect to}
\def\Wlog/{Without loss of generality}
\def\a{\alpha}

\def\Ab/{Abel}
\def\Ban/{Banach}
\def\Bansp/{\Ban/ space}
\def\Belt/{Bel\-tra\-mi}
\def\Berg/{Berg\-man}
\def\Bern/{Ber\-nou\-lli}
\def\Berz/{Berezin}
\def\Bess/{Bessel}
\def\Cart/{Car\-tan}
\def\Cay/{Cay\-ley}
\def\CG/{Clebsch-Gordan}
\def\Cl/{Clifford}
\def\CR/{Cauchy-Rie\-mann}
\def\Dir/{Dirichlet}
\def\Eucl/{Euclide}
\def\Eucln/{Euclidean}
\def\F/{Fourier}
\def\Hank/{Hankel}
\def\Hankf/{\Hank/ form}
\def\Herm/{Hermite}
\def\Hilb/{Hilbert}
\def\Hilbs/{Hilbert space}
\def\Hilbsp/{Hilbert space}
\def\HS/{Hilbert-Schmidt}
\def\Lag/{La\-grange}
\def\Lap/{La\-place}
\def\LapBelt/{\Lap/-\Belt/}
\def\Leb/{Lebesgue}
\def\Marc/{Mar\-cin\-kie\-wicz}
\def\Moeb/{Moebius}
\def\Moebt/{Moebius transformation}
\def\Moebtransfn/{Moebius transformation}
\def\Pla/{Plan\-che\-rel}
\def\Poin/{Poin\-car\'e}
\def\Riem/{Rie\-mann}
\def\Riemn/{\Riem/ian}
\def\psRiemn/{pseudo-\Riem/ian}
\def\Riems/{Rie\-mann surface}
\def\Schroe/{Schr\"odinger}
\def\Weier/{Weier\-strass}
\def\anal/{analytic}
\def\bsd/{bounded symmetric domain  }
\def\bdd/{bounded}
\def\calc/{calculation}\def\conj{conjugate}
\def\calci/{calculating}\def\eg{e.g.}
\def\conj/{conjugate}
\def\deco/{decomposition}
\def\eg/{e.g.}
\def\fct/{function}
\def\gp/{group}
\def\hw/{highest weight}
\def\hwv/{highest weight vector}
\def\hwvs/{highest weight vectors}
\def\lw/{lowest weight}
\def\lwv/{lowest weight vector}
\def\lwvs/{lowest weight vectors}
\def\hds/{holomorphic discrete series}
\def\iff/{if and only if}
\def\inv/{invariant}
\def\irrde/{irreducible decomposition}
\def\meas/{measure}
\def\transf/{transform}
\def\rep/{representation}
\def\resp/{respectively}
\def\inters/{intertwines}
\def\interg/{intertwining}
\def\meta/{metaplectic representation}
\def\qu/{quaternion}
\def\rep/{representation}
\def\symdom/{ symmetric domain}
\def\st/{such that}
\def\shd/{subhead}
\def\transf/{transform}
\def\wrt/{with respect to}

\def\Norm#1#2#3{\Vert#1\Vert^{#3}_{{#2}}}

\begin{document}
\def\abstractname{Abstract}
\def\chrefname{References}

\title[Tensor products of complementary series]{Tensor 
products of complementary series
of rank one Lie groups
 }
\author{ Genkai Zhang}

\address{Mathematical Sciences, Chalmers University of Technology and
Mathematical Sciences, G\"oteborg University, SE-412 96 G\"oteborg, Sweden}
\email{\it genkai@chalmers.se}
\thanks{This research was partially supported by 
the Swedish
Science Council (VR)}

\begin{abstract}
We consider the tensor product
$\pi_{\alpha}\otimes \pi_{\beta}$ of complementary
series  representations $\pi_{\alpha}$ and 
$\pi_{\beta}$  of classical rank one groups
$SO_0(n, 1)$,
$SU(n, 1)$ and 
$Sp(n, 1)$.
We prove that there is a discrete 
component $\pi_{\alpha+\beta}$  for small
parameters $\alpha, \beta$ (in our parametrization). 
We prove further that for $SO_o(n, 1)$
there are  finitely many
complementary series of the form $\pi_{\alpha+\beta + 2j}$, $j=0, 1,
\cdots, k$, 
appearing in
the tensor product
$\pi_{\alpha}
\otimes \pi_{\beta}
$ of  two
complementary series $\pi_{\alpha}$
and $\pi_{\beta}$, where
$k=k(\alpha, \beta, n)$ depends on $\alpha, \beta, n$. 
\end{abstract}

\maketitle

\baselineskip 1.35pc

\section{Introduction}

The question of finding discrete components 
of
a unitary representation $(\tau, L)$ of a Lie group $L$
under a Lie subgroup $G\subset L$ is one of the central topics in representation theory. 
As a problem in harmonic analysis
this would involve some explicit
realizations
of
 representations for $L$ and $G$, 
the existence and boundedness of
intertwining operators from 
$(\tau, L)$
to  discrete components $(\pi, G)$.
 The problem
of finding {\it formal} intertwining operators
from the smooth representations
$(\tau^\infty, L)$
to $(\pi^\infty, G)$ can be studied
somewhat by algebraic methods. Given the
intertwining operators  the next question
is to prove their boundedness. 
In the present paper we shall consider
the case when $(\pi_{\gamma}, G)$ 
is a complementary
series of a 
classical rank one group $G$ and $(\tau, L)=
(\pi_{\alpha}\otimes \pi_{\beta}, G\times G)$
is a tensor product of the same type. We 
  prove the existence
of discrete components 
$(\pi_{\gamma}, G)$  in  the tensor product.

The complementary series are unitarization
of principal series and can be realized as function
spaces on the Euclidean space $\mathbb R^{n-1}$
or $H$-type groups.
There exist
certain  bilinear invariant differential
intertwining operators on the tensor products, and the study
of them is of interesting for various purposes.
 The most
well-known examples of similar operators might
be  the Rankin-Cohen brackets
on  tensor products
of holomorphic discrete series of $\SL(2, \mathbb R)$, which yield
also a decomposition of the tensor products
in the unitary sense. 
In \cite{Ovsienko-Redou}  Ovsienko and Redou
have
found a family of  formal differential operators
on tensor product of  principal series
$\pi_{\alpha}, \alpha\in \mathbb C$,
of the conformal group $SO(n, 1)$.
We shall prove the boundedness for the intertwining
operators under certain conditions, thus proving
the existence of finitely many discrete components. For
the group $SU(n, 1)$ and $Sp(n, 1)$
we prove the natural restriction operator
is bounded and found then one discrete component, namely
$\pi_{\alpha+\beta}$ in 
$\pi_{\alpha} \otimes\pi_{\beta}$. Our main results 
are stated in Theorems 4.1 and 5.1.

We give a very brief account of the closely related
results. 
The  Rankin-Cohen brackets for holomorphic
representations  and their generalizations
have been studied in various contexts. An analytic
and effective method to derive
them is via the reproducing formula. 
The similar tools for spherical representations
are the Knapp-Stein intertwining operators.
The intertwining operators in \cite{Ovsienko-Redou} are
 found  by using an Ansatz
expressing the operators as 
polynomials of the Laplacian operators $\mathcal L_x$, $\mathcal L_y$,
and the inner product $\nabla_x\cdot \nabla_y$ evaluated 
on the diagonal $x=y$.
The same operators are obtained in \cite{Clerc-Beckmann}
as residues of a family of integral 
bilinear intertwining operators, which
is based on some earlier works
on trilinear form \cite{Clerc-2011,
Clerc-Orsted, CKOP, BKZ}, which are variations
of the Knapp-Stein intertwining operator.
Kobayashi and Pevzner \cite{Kob-Pev} have
constructed systematically
invariant bilinear
differential operators using special orthogonal polynomials.
Tensor product
decompositions for representations of $SL(2, \mathbb R)$
had been studied in greater details
in the work of Pukanszky
\cite{Pukanszky62} and by Repka
\cite{Repka78};  see  also \cite{Asmuth-Repka}.
Theses
results combined with the
general theory of Burger-Li-Sarnak \cite{Burger-Li-Sarnak}
have also found applications in automorphic forms \cite{Clozel-ias-park-lec}.
Tensor products of representations of $SL(2, \mathbb C)$ (i.e.,
locally isomorphic to the Lorentz group $SO_0(3, 1)$) have
been studied by Naimark \cite{Naimark-1961}; see also
\cite{Neretin-1995} where some complementary series
representations were constructed using  restriction
of holomorphic representations, the same
idea being used in the present paper in the construction
of discrete components for the groups
$SU(n, 1)$ and $Sp(n, 1)$. A related question
is the  branching of complementary
series of rank one groups $L=SO(n, 1)$,
$SU(n, 1)$, $Sp(n, 1)$ under 
the subgroups $G=SO(n-1, 1)$,
$SU(n, 1)$, $Sp(n, 1)$ and is studied in  \cite{Mollers-Oshima, Speh-Zhang, Zhang-jfa}.

The paper is organized as follows.
In Section 2 we recall the realization
of the complementary series of $SO(n, 1)$ as function spaces on 
$\mathbb R^{n-1}$. For the groups
$G=SU(n, 1)$ and $Sp(n, 1)$ we shall
need their realization in terms of holomorphic
representations of some  Hermitian groups containing $G$.
In Section 3 we present an elementary 
and independent construction
of the bilinear differential operators for
$SO(n, 1)$. While the technical
computations are slightly different from those
in the existing literature, the underlying idea
is still the same, the explicit application of Knapp-Stein 
operators here making the computations conceptually clear.
Section 4 is devoted to the proof of
the appearance of finitely many discrete components,
$\pi_{\alpha + \beta+2j}$, $j\ge 0$,
in the tensor product $\pi_{\alpha}\otimes
\pi_{\beta}, \alpha, \beta>0$,
 when the parameters $\alpha$ and $\beta$ 
are relatively small (in our parametrization);
In \S5 we treat the  other rank one groups $SU(n, 1)$
and $Sp(n, 1)$.

It is clear now that there remain
many interesting problems in understanding
the discrete components of the tensor products
of  complementary series, in particular
for the other rank one groups $SU(n, 1)$,
$Sp(n, 1)$ and exceptional $F_4$-group.
The spherical representations
are then realized on H-type groups
and the Euclidean Fourier analysis in Section 4 might have to be
replaced by Fourier analysis on the nilpotent groups.
First of all it would be interesting
to construct bilinear invariant
differential operators.
Secondly one may proceed to prove
the boundedness and it seems
estimating these operators becomes far more complicated.
On the hand one might also try
the method of holomorphic
extension as in the present paper.

I would like to thank Jean-Louis Clerc 
 for some stimulating discussions.

\section{Spherical representations of rank one group $G$}
We fix notation and 
 recall  some known results on
induced representations of $G$ 
and the Knapp-Stein intertwining operator.
We shall use the  non-compact realization
of the representations. 
We shall be quite brief, and most of the technical
 formulas can be found e.g. in 
\cite{Johnson-Wallach,CDKR-jga}
where  the general case of rank one groups
is studied.

\subsection{Classical rank one Lie groups}
Let $G=SO_0(n, 1, \mathbb F)=SO_0(n, 1), SU(n, 1), Sp(n, 1)$,
 $\mathbb F=\mathbb R, \mathbb C, 
\mathbb H$, 
be the classical rank one Lie group
in its standard realization \cite{Johnson-Wallach,
Speh-Zhang}.
Let $\fg=\fk +\fp$ be the Cartan decomposition of the Lie algebra
$\fg$
of $G$.
We fix an element $H\in \fp$
and the subspace $\fa:=\mathbb RH_0\in \fp$
such that $\Ad(H)$ has (possible) eigenvalues $\pm 2,
\pm 1, 0$.
The root space decomposition of $\fg$ under  $H$ 
is 
$$
\fg=\fn_{-2} + \fn_{-1} +
(\fa+ \fm) + \fn_{1} +\fn_{2} 
$$
with eigenvalues $\pm 2, \pm 1, 0$
if $\mathbb F=\mathbb C, \bh$,
and $\pm 1, 0$ if $\mathbb F=\mathbb R$.
  Here $\fm\subset \fl$ is the zero root space.
We denote by $\fn=\fn_1 \oplus \fn_2$ the sum of the positive root spaces.
 Then
$\fm+\fa +\fn$
is a maximal parabolic subalgebra
of $\fg$.

Let $\rho$
be the half sum of positive roots. Then
$$
\rho(H)=
\begin{cases}
\frac{n-1}2, &\quad \mathbb F=\br\\
n, &\quad \mathbb F=\bc\\
2n+1, &\quad \mathbb F=\bh
\end{cases}
$$
and we write  $\rho=\rho(H)$.

\subsection{Spherical representations and complementary series
for $G=SO_0(n, 1; \mathbb F)$}

Denote $M, A, N$ 
the corresponding 
subgroups of $G$ with Lie algebras
$\fm, \fa, \fn$, and $P=MAN$ the parabolic subgroup.
For $\mu\in \mathbb C$ 
let $\pi_{\mu}^\infty$
 be the 
induced {\it smooth}
representation
of $G$ from the character $
e^{-\mu}:
me^{tH}
n \in P=MAN\mapsto e^{-\mu t}$, consisting of $C^\infty$-functions $f$ on $G$
such that
\begin{equation}
  \label{eq:ind-norm}
f(g 
me^{tH_0}
n
)=
e^{-\mu t}
f(g),  me^{tH_0}n\in MAN.  
\end{equation} 
In particular  $f$ are determined by their
restriction on $K$ and are identified further
as smooth functions on  $K/M$. We have  $\pi_{\mu}^\infty
=C^\infty(K/M)=C^\infty(S)$ as vector spaces. Restricting 
smooth functions in 
$\pi_{\mu}^\infty$
 to $N^-$ results in an injective
map to a subspace of $C^{\infty}(N^-)
=C^{\infty}(\mathbb R^{n-1})$. We shall fix this realization
of $\pi_{\mu}^\infty$. (Indeed it is only a proper subspace
of $C^{\infty}(N^-)$ as some matching conditions at infinity are needed.)

The explicit formulas for $\pi_\mu(g)$ can
be found in \cite{Johnson-Wallach} in the compact
picture and in
\cite{Speh-Venk-2} for the non-compact picture. We shall only need
them for the real group
$G=SO_0(1, n)$.
In this case $\mathbb n=\mathbb R^{n-1}$ and $G$ is generated by
the parabolic group
 $ MAN^-$ 
and 
the Weyl group element
$w$, \cite[Theorem 1.4, Ch.~IX]{He2},
where  $ MAN^-$ acts on $\mathbb R^{n-1}$
by the defining action
and $w$ by inversion, $w(x)=-\frac{1}{x}:=
-\frac{x}{|x|^2}$.
Their actions on $\pi_{\mu}^\infty$
 are  given by
$$
\pi_{\mu}(g)f(x)= e^{-t\mu} f(e^{t}m^{-1}(x-x_0)),
\quad (m, e^{tH}, x_0)
=me^{tH} x_0
\in MAN^-, \quad N^{-}=\mathbb R^{n-1},
$$
and
$$
\pi_{\mu}(w)f(x)= \Vert x\Vert^{-2\mu }
f(-\frac{x}{\Vert x\Vert^2}).
$$
Note also that the Jacobians of 
$g=(m, e^{tH}, x_0)$ and  of the Weyl group element $w$ 
on $N^-=\mathbb R^{n-1}$ 
are given by
\begin{equation}
\label{eq:jac}
J_g(x)=e^{t(n-1)},
\quad
J_w(x)=\frac{1}{|x|^{2(n-1)}}.
\end{equation}

The representation  $\pi_{\mu}(g), g\in G$,
$\mu\in (\rho +i\mathbb R)$ is already unitary
for the natural unitary norm in $L^2(K/M)$. However
for certain real $\mu\in \mathbb R$ 
it is possible to define
a different $\fg$-invariant inner product on the space of $K$-finite
vectors, in which case its Hilbert space
completion is called {\it complementary series}; see
\cite{Johnson-Wallach}. The precise range of the parameters $\mu$ is given by
\begin{equation} \label{mu-fixed}
\mu\in  \begin{cases} 
(0, 2\rho), \mathbb F=\br, \bc
\\
(2, 2\rho-2), \mathbb F=\bh.
\end{cases} 
\end{equation} 
 We shall
denote the corresponding representation
still by $\pi_{\mu}$, and shall
 use its non-compact realization for the real case $G=SO_0(n,
 1)$, allowing us to find  (generically) more than one
 discrete components in the tensor product decomposition in \S4.

\subsection{Realization
of complementary series
for $G=SO_0(n, 1)$ on $\mathbb R^{n-1}$}

The unitarization of the complementary series
 is obtained via the Knapp-Stein intertwining operator,
defined preliminarily on $K$-finite vectors (which
can be obtained from $K$-finite vectors on $K/M$ via Cayley transform),
\begin{equation}
  \label{eq:j-n-c}
J_{\mu} f(x)
=\int_{\mathbb R^{n-1}} K_{\mu}(x, y) 
f(y) dy,  
\end{equation}
where 
\begin{equation}
  \label{eq:j-n-c-k}
 K_{\mu}(x, y) = C_{\mu}
\frac1 {|x-y|^{2\mu}}, \quad
C_\mu
=\frac
{\Ga(\rho-\frac{\mu}2)
\Ga(\rho-\frac{\mu}2 +\frac 12)}
{\Ga(\frac{n}2)
\Ga(\rho-\mu)}.
\end{equation}
(The normalization is chosen here so that in the compact picture
$J_\mu 1_S=1_S,$ where $1_S$ is the constant function on $S$ viewed
as a function on $G$ restricted to $\mathbb R^{n-1}.$)

Then $J_\mu$ is a $G$-intertwining operator
$$
\boxed{J_\mu: \pi_{\widetilde \mu}^\infty
\to \pi_{\mu}^\infty, \,
 \quad \widetilde \mu :={2\rho-\mu}}
$$
for $\mu <<0$.  
It has holomorphic continuation
to the whole complex plane, and in particular  holomorphic and non-zero in the two symmetric
strips around $\Re \mu=\rho $,
\begin{equation}
  \label{eq:hol-str}
\{\mu; 0<\Re\mu < \rho\},
\quad
\{\mu; \rho <\Re\mu < 2\rho
\}.
\end{equation}
The formal intertwining property
can be proved by using the following transformation
rule of $K_\mu$,
$$
K_\mu(gz, gw)=  |(cz +d)|^{-\mu}
K_\mu(z, w)|(cw +d)|^{-\mu} 
=  J_g(z)^{\frac{\mu}{n-1}} 
K_\mu(z, w)
 J_g(w)^{\frac{\mu}{n-1}} 
$$
where $J_g$ is the Jacobian of the action of $g\in G$ on
$N^-=\mathbb R^{n-1}$. See \cite[Chapter VII]{Kn-book}
and \cite{Vogan-Wallach}
for some general
theory of intertwining operators.

 The inner product
\begin{equation}
  \label{eq:j-n-c-norm}
(f_1, f_2)_\mu=(J_{\widetilde \mu}f_1, f_2)_{L^2(\mathbb R^{n-1})}
\end{equation}
for $f_1, f_2\in  C^\infty_0(\mathbb R^{n-1})$
is  a pre-Hilbert norm, and is invariant under $g\in G$ sufficiently
close
to the identity (depending on 
$f_1, f_2$). The completion  defines the
complementary series, $\mu\in (0, 2\rho)$.
We shall use its description 
in term of Fourier transform $f\mapsto \mathcal Ff$.
The space  $\pi_\mu$  is the completion
of $C_0^\infty(\mathbb R^{n-1})$
under  the (equivalent) norm
\begin{equation}
  \label{eq:Four-norm}
\Vert 
f\Vert_{\mu}^2
=\int_{\mathbb R^{n-1}}
|\mathcal F f(\xi)|^2 |\xi|^{n-1-2\mu} d\xi
=\Vert 
\mathcal F f(\cdot) 
|\cdot |^{\frac{\widetilde\mu}{2}} 
\Vert^2_{
L^2(
\mathbb R^{n-1}
)
},
\end{equation}
for $0<\mu <2\rho$.
See e.g. \cite{Speh-Venk-2}.

\subsection{Complementary series
for $G=SU(n, 1), Sp(n, 1)
$
and their holomorphic extensions}

We shall use a different method for  the cases
 $G=SU(n, 1), Sp(n, 1)$. The method is based on, roughly
speaking, {\it holomorphic extension}.

We consider a Hermitian Lie group
$G_1$ containing $G$
as a symmetric subgroup such that
$G/K$ is a real form of the
Hermitian subgroup $G_1/K_1$; see \cite{Dijk-Hille-jfa, Speh-Zhang}.
In certain rather general context
the similar groups $G_1$ are also called {\it overgroups}
 in some Russian litterature; see e. g. \cite{Molchanov}.

More precisely let 
$$
G_1=\begin{cases}
SU(n, 1)\times SU(n, 1), \mathbb F=\mathbb C\\
SU(2n, 2), 
 \mathbb F=\bh,
\end{cases}
$$
with $G=SU(n, 1), Sp(n, 1)$ being
realized as the
diagonal  subgroup of $G_1$
and respectively as complex 
transformations via the standard identification
of $\mathbb H=\mathbb C^2$.

The holomorphic discrete series
of $G_1$ can be realized on
the space of holomorphic functions
on $D_1$. To fix notation we let
$$
V_1=\begin{cases} 
\bc^n \oplus \overline{\bc^n},  & \mathbb F=\mathbb C\\
\\
 M_{2n, 2}(\mathbb C), & \mathbb F=\bh,
\end{cases}
$$
and the space 
$D_1=G_1/K_1$
is realized as a bounded symmetric domain in $V_1$,
$$
D_1=\begin{cases} 
B^n \times \overline{B^n}, &  \mathbb F=\mathbb C
\\
\{Z\in M_{2n, 2}(\mathbb C); Z^*Z<I\}, & 
 \mathbb F=\bh.
\end{cases}
$$
 Let 
$\mathcal H_\nu(D_1)$ 
be the space
of holomorphic functions on $D_1$ with reproducing kernel
$h(z, w)^{-\nu}$ for $\nu$ sufficiently large,
where $$
h(z, w)=
\begin{cases} 
(1-\langle z_1, w_1\rangle )
(1-\langle w_2, z_2\rangle ),
&  \mathbb F=\mathbb C
\\
\det(1- w^\ast z), &  \mathbb F=\mathbb H.
\end{cases}
$$

It is now well-known  \cite{FK} that if $\nu$ is in the set
$$
\begin{cases} (0, \infty) & \mathbb F
=\mathbb R, 
 \mathbb C \\
(1, \infty) &  \mathbb F= \mathbb H,
\end{cases}
$$
then the kernel $h(z,w)^{-\nu}$
is positive definite   and it defines
a unitary (projective) representation
of $G_1$.
We denote this representation by
$(\mathcal H_\nu(D_1), \tau_\nu, G_1)$. 

We recall the following theorem
\cite{Dijk-Hille-jfa}.
\begin{theo+} 
The complementary series
$(\pi_{\mu}, G)$ appears as 
a discrete summand 
in  
 $(\mathcal H_\nu(D_1),  \tau_\nu, G_1)$ restricted to $G$ if $\nu$ and $\mu$ are
related by
$$2\nu 
=
 \begin{cases} 
\mu, \,\,  \mu \in (0,  n), \quad  & \mathbb F=\bc \\
\mu, \,\, \mu\in (2, {2n-1}), \quad & \mathbb F=\bh.
 \end{cases}
$$
\end{theo+}

Note that the range of $\mu$ is, disregarding the Weyl group symmetry,
precisely the whole range of the complementary series representations.
In other words, any complementary series of $G$ is a discrete component
in an holomorphic representation of $G_1$.

\section{Invariant  bilinear differential operators for general
spherical series representations of $G=SO_0(n, 1)$}

We denote $\pi_{\alpha}^\infty
 \otimes
\pi_{\beta}^\infty$ the induced smooth representation 
of $G\times G$ from the parabolic subgroup 
$P\times P$; see (\ref{eq:ind-norm}). (The tensor notation here
is a bit ambiguous, and the precise definition should be that it
is the induced
representation from the tensor product of representations of $P$.)
The group $G$ is viewed as the diagonal subgroup
of $G\times G$.

\begin{theo+} For any nonnegative integer $m\ge 0$
there exists
 a    $G$-intertwining  differential operator 
$\mathcal D_{\alpha, \beta, m}
$ of degree $2m$
 meromorphic  in $(\alpha, \beta)\in \mathbb C^2$,
$$
\mathcal D_{\alpha, \beta, m} :
\pi_{\alpha}^\infty
 \otimes
\pi_{\beta}^\infty
 \to  \pi_{\alpha+\beta +2m}^\infty.
$$
The  only possible
poles of $
\mathcal D_{\alpha, \beta, m}
$ appear
when  $\alpha$ or $\beta\in \Lambda_m$, where 
$$
\Lambda_j=\{0, -1, \cdots, -m+1\} \cup 
\left( \rho -1  +\{0, -1, \cdots, -m+2\}
\right).
$$
\end{theo+}

The proof will be divided into  a few elementary steps.

Let $S_{\alpha, \beta, m}
(x, y; z, w) $ be the kernel
\begin{equation}
  \label{eq:S-ker}
S_{\alpha, \beta, m}
(x, y; z, w)
=
\left(
    \frac{|(x-z)-(y-w)|^2}
{|x-z|^{2}
|y-w|^{2}
}
\right)^m
  \frac{1}
{|x-z|^{2\alpha}
|y-w|^{2\beta}
},
\end{equation}
and write for simplicity 
$$
S_{\alpha, \beta, m}(x; z, w)=S_{\alpha, \beta, m}(x, x; z, w),
$$
the restriction  to the diagonal $x=y$  of 
$S_{\alpha, \beta, m}(x, x; z, w)$.
\begin{lemm+}  The integral operator
$$
T_{m}f(x)=T_{\alpha, \beta, m}f(x):
=
C_\alpha
C_\beta
\int_{\mathbb R^{2(n-1)}}
S_{\alpha, \beta, m}(x; z, w)
 f(z, w) dz dw
$$
defines  an intertwining operator 
$$
\pi_{\widetilde \alpha}^\infty
\otimes \pi_{\widetilde \beta}^\infty
\to \pi_{\alpha+\beta+2m}^\infty,
$$
for $\alpha, \beta\in \mathbb R$ sufficiently negative
and has meromorphic continuation
to all 
$\alpha, \beta\in \mathbb C$.
\end{lemm+}
\begin{proof}
Recall the  group $G$ is generated by $P$
and $w$.
The formal intertwining property
follows directly from
a  change of variables $(x, y)\mapsto (gx, gy)$
for $g\in P$ and $g=w$
along with the formula 
(\ref{eq:jac})
for the Jacobians.
To prove the meromorphic 
continuation in $\alpha$ and $\beta$
we observe that changing
$(x, y)$ to $(x-z, y-z)$ we need only
to prove that the integral
$$
\int_{\mathbb R^{2(n-1)}}
   \frac{|x-y|^{2m}
}
{|x|^{2m}
|y|^{2m}
}
  \frac{1}
{|x|^{2\alpha}
|y|^{2\beta}
}
 f(x, y) 
dx dy
$$ 
is 
meromorphic in $(\alpha, \beta)$. But this
is just up to normalization  constants
the integral $(J_{\alpha +m} \otimes 
J_{\beta +m})(F)$, $F(x, y)=
|x-y|^{2m} f(x, y) $ and thus has the
continuation.
\end{proof}

In the compact-realization this operator is
$$
T_m f(x)=C_\alpha 
C_\beta 
\int_{S\times S}
    \frac{
(1-\langle z, w\rangle)^m
}
{
(1-\langle x, z\rangle)^{m+\alpha}
(1-\langle x, w\rangle)^{m+\beta}
}
f(z, w)dz dw.$$
That the integral is well-defined for $\alpha, \beta <<0$
 can also be  deduced from this formula.

Next  we need some elementary known
 identities
for the Laplacian operator $\mathcal L=\partial_1^2 
+\cdots +\partial_n^2$ acting on $|x|^{-2\alpha}$,
the so-called  Bernstein-Sato identities.
Recall   the Pochammer symbol defined by $(\alpha)_j=\alpha
(\alpha+1)\cdots (\alpha +j-1)$.
\begin{lemm+}
 The following differentiation  formula holds 
\begin{equation}
  \label{eq:L-on-K}
\mathcal L^j |x|^{-2\alpha}
=2^{2j}(\alpha)_j
(\alpha+1-\rho)_j
|x|^{-2(\alpha+j)}, \quad x\ne 0.
\end{equation}
\end{lemm+}

We define inductively a family of differential
operators of constant coefficients on $C^\infty(\mathbb R^{2(n-1)})$
by $$
M_{\alpha, \beta, 0}=I,\quad  M_{\alpha, \beta, 1}=\nabla_x \cdot
\nabla_y$$
 and
 \begin{equation}
   \label{eq:D-N-M}
\mathcal M_{\alpha, \beta, j+1}
=(\nabla_x \cdot \nabla_y) 
\mathcal M_{\alpha, \beta, j}
-\frac {j( n-1-3j-2 \alpha -2\beta )
 } 
{
 (\alpha +1-\rho)
 (\beta +1-\rho)
}
\mathcal M_{ \alpha+1, \beta+1, j-1
}
\mathcal L_x
\mathcal L_y.   
 \end{equation}
It follows from the construction that the only possible poles of 
$M_{\alpha, \beta, j}$, $j\ge 2$,  appear
when   $\alpha $ or     $\beta$ is in 
$$
\{\rho -i; i=1, \cdots, j-1\}.
$$

\begin{lemm+} 
The following formula holds for all $(\alpha, \beta)\in \mathbb C^2$
and $m=0, 1, 2, \cdots$,
\begin{equation}
  \label{eq:M-on-kernel}
  \mathcal M_{\alpha, \beta, m}
\left(\frac{1}{|x|^{2\alpha}|y|^{2\beta}}\right)
= 2^{2j}(\alpha)_m (\beta)_m
\left(
 \frac{\langle x, y\rangle
}
{
|x|^{2}
|y|^{2}
}
\right)^m
\left(\frac{1}{|x|^{2\alpha}|y|^{2\beta}}\right).
\end{equation}
\end{lemm+}
\begin{proof} We write
the identify
as $\text{LHS}_m=\text{RHS}_m$
prove it  by induction on $m$.

The identity  is trivially true for $m=0$.
Assuming the identity holds for  $0\le m\le j$
for all $\alpha, \beta$.
We perform the differentiation.
$\nabla_x \cdot \nabla_y$ on the identity with $m=j$ and find
\begin{equation}
  \label{eq:induct}
\nabla_x \cdot \nabla_y  \mathcal M_{\alpha, \beta, j}
\frac{1}{|x|^{2\alpha}|y|^{2\beta}}
=2^{2j}(\alpha)_j (\beta)_j
(I +II)=
2^{2j}(\alpha)_j (\beta)_j 
I
+2^{2j}(\alpha)_j (\beta)_jII
\end{equation}
a sum of two terms, with the first term 
\begin{equation*}
\begin{split}
2^{2j}(\alpha)_j (\beta)_j I&=2^{2j}(\alpha)_j (\beta)_j
2^2(\alpha +j)(\beta +j)
\left(
 \frac{
\langle x, y\rangle
}
{
|x|^{2}
|y|^{2}
}
\right)^{j+1}
\frac{1}{|x|^{2\alpha} |y|^{2\beta}
}
\\
&=
2^{2(j+1)}(\alpha)_{j+1} (\beta)_{j+1}  
\left(
 \frac{
\langle x, y\rangle
}
{
|x|^{2}
|y|^{2}
}
\right)^{j+1}
\left(\frac{1}
{
|x|^{2\alpha} |y|^{2\beta}
}\right)
=\text{RHS}_{j+1},
\end{split}
\end{equation*}
the the RHS of  (\ref{eq:M-on-kernel}) for $m=j+1$,
and 
$$
II=
j (n-1-3j-2\alpha -2\beta) 
\left(
 \frac{
\langle x, y\rangle
}
{
|x|^{2}
|y|^{2}
}
\right)^{j-1}
\left(
\frac{1}
{
|x|^{2(\alpha+1)} 
|y|^{2(\beta+1)}
}
\right).
$$
We treat the second term using the induction hypothesis for $m=j-1$
with $(\alpha,\beta)
$ being replaced
by  $(\alpha+1,\beta+1)$,
\begin{equation*}
  \begin{split}
&\quad 2^{2(j-1)}
(\alpha+1) _{j-1}
(\beta+1) _{j-1}
\left(
 \frac{
\langle x, y\rangle
}
{
|x|^{2}
|y|^{2}
}
\right)^{j-1}
\left(
\frac{1}{|x|^{2(\alpha+1)} |y|^{2(\beta+1)}}
\right)
\\
&=\mathcal M_{\alpha+1,\beta+1,
j-1 
}
\left(
\frac{1}{
|x|^{2(\alpha+1)} |y|^{2(\beta+1)}}
\right)
  \end{split}
\end{equation*}
which is furthermore
$$
\frac 1{2^2\alpha (\alpha +\rho-1)
\beta (\beta +\rho-1)
}
\mathcal M_{
\alpha+1,\beta+1,
j-1}
\mathcal L_x
\mathcal L_y 
\left(
\frac{1}{|x|^{2(\alpha)} |y|^{2(\beta)}}
\right).
$$
We rewrite 
(\ref{eq:induct})  as
\begin{equation*}
  \begin{split}
\text{RHS}_{j+1}&=(\nabla_x \cdot \nabla_y  \mathcal M_{
\alpha, \beta, j}
-\frac {j(n-1-3j-2\alpha -2\beta) 
}
{ (\alpha +\rho-1)
 (\beta +\rho-1)
}
\mathcal M_{
\alpha+1,\beta+1,
j-1}
\mathcal L_x
\mathcal L_y)
\frac{1}{|x|^{2(\alpha)} |y|^{2(\beta)}}\\
&=\mathcal M_{
\alpha, \beta, j+1 
}
\left(
\frac{1}{|x|^{2(\alpha)} |y|^{2(\beta)}}
\right)\\
&=\text{LHS}_{j+1}
  \end{split}
\end{equation*}
by the definition. This finishes the proof.
\end{proof}

Combining the two Lemmas  we have
\begin{equation*}
\mathcal M_{\alpha+j, \beta+i, k}
\mathcal L_x^j
\mathcal L_x^i
   \frac{1}
{|x|^{2\alpha}
|y|^{2\beta}
}
=c_{i, j, k}(\alpha, \beta)   
\left(
 \frac{\langle x, y\rangle}
{
|x|^{2}
|y|^{2}
}
\right)^k
   \frac{1}
{|x|^{2\alpha+2j}
|y|^{2\beta+2i}
},
\end{equation*}
where 
\begin{equation}
c_{i, j, k}(\alpha, \beta)
=  2^{2k+2j+2i} 
(\alpha)_{j+k}
(\alpha+1-\rho)_j
(\beta)_{i+k}
(\beta+1-\rho)_i.
\end{equation}
Here we have used the fact that
$$
(\gamma)_{j} 
(\gamma+j)_{k}
=(\gamma)_{j+k}.
$$
By translation invariance
we have 
\begin{equation}
  \label{eq:M-L-on-kernel}
\begin{split}
&\quad\,\mathcal M_{
\alpha+j, \beta+i, 
k}
\mathcal L_x^j
\mathcal L_x^i
   \frac{1}
{|x-z|^{2\alpha}
|y-w|^{2\beta}
}
\\
&=c_{i, j, k}(\alpha, \beta)   
\left(
 \frac{\langle x-z, y-w\rangle}
{
|x-z|^{2}
|y-w|^{2}
}
\right)^k
   \frac{1}
{|x-z|^{2\alpha+2j}
|y-w|^{2\beta+2i}
},
\end{split}
\end{equation}

We prove now Theorem 3.1.
\begin{proof} 
The operator 
$$
T_{\alpha, \beta, m} 
(J_{\alpha}\otimes  
J_{\beta}):
\pi_{\widetilde \alpha}^\infty
\otimes \pi_{\widetilde \beta}^\infty
\to \pi_{\alpha+\beta+2m}^\infty
$$
 is an intertwining operator by Lemma 3.2. We prove
it is  a differential operator.
The idea is
to  differentiate 
 the identity $f=
(J_\alpha\otimes 
J_\beta)
(J_{\widetilde \alpha}\otimes 
J_{\widetilde \beta}) f$. 
We shall perform formal computations
on the integral first and justify the procedure
in the end. 
Let $f\in \pi_{\alpha}^\infty\otimes
\pi_{\beta}^\infty$ and $
g=
(J_{\widetilde \alpha}\otimes
J_{\widetilde \beta} )f.
$
We denote
\begin{equation*}
\begin{split}
  \label{eq:2}
&\quad\, \mathcal E_{\alpha, \beta, m}
f(x, y)\\
&=
\sum_{i+j +k=m} 
\varepsi_{i, j, k}(\alpha, \beta)
\mathcal 
M_{\alpha +j,
\beta +i, k}
\mathcal L_x^{j}
\mathcal L_y^{i} f(x, y)
\end{split}
\end{equation*}
and
\begin{equation}
  \label{eq:D-N}
\mathcal D_{\alpha, \beta, m}f(x)=
\mathcal E_{\alpha, \beta, m}f{|}_{x=y},
\end{equation}
for $f\in C^\infty(\mathbb R^{2(n-1)})$,
where
$$
\varepsi_{i, j, k}(\alpha, \beta)
:=\binom{m}{i, j, k}
\frac{(-2)^{k}}
{c_{i, j, k}(\alpha, \beta)
}.
$$
We claim
that 
\begin{equation}
  \label{eq:D-N-TJ}
\mathcal D_{\alpha, \beta, m} 
f=
T_{\alpha, \beta, m} 
(J_{\alpha}\otimes  
J_{\beta})f, \quad f\in \pi_{\widetilde \alpha}^\infty
\otimes \pi_{\widetilde \beta}^\infty
\end{equation}
proving the formal intertwining property of the differential
operator $ \mathcal D_{\alpha, \beta, m} $.

The binomial expansion of 
$S(x, y; z, w)$
reads as follows
\begin{equation*}
  \begin{split}
S(x, y; z, w)
&=\left(\frac{
|x-z|^2
+|y-w|^2
  -2
\langle x-z, y-w\rangle}
{|x-z|^2 |y-w|^2
}
\right)^m
  \frac{1}
{|x-z|^{2\alpha}
|y-w|^{2\beta}
}
\\
&=\sum_{i+j +k=m}
\binom{m}{i, j, k}
{(-2)^{k}}
\left(\frac
{\langle x-z, y-w\rangle}
{
|x-z|^{2}
|y-w|^{2}
}
\right)^{k}
  \frac{1}
{
|x-z|^{2j+2\alpha}
 |y-w|^{2i+2\beta} 
}.
  \end{split}
\end{equation*}
Summing the formula 
(\ref{eq:M-L-on-kernel})
over $(i, j, k)$ we have then
$$
\mathcal E_{\alpha, \beta, m}
   \frac{1}
{|x-z|^{2\alpha}
|y-w|^{2\beta}
}
=S_{\alpha, \beta, m}
(x, y; z, w)
$$
which further implies that
\begin{equation}
  \label{eq:d-n-a-b}
\mathcal D_{\alpha, \beta, m}
\left(   \frac{1}
{|x-z|^{2\alpha}
|y-w|^{2\beta}
}\right)
=S_{\alpha, \beta, m}(x, x; z, w) =S_{\alpha, \beta, m}(x; z, w).
\end{equation}

The identity $f=(J_{ \alpha}\otimes
J_{ \beta})
(J_{\widetilde \alpha}\otimes
J_{\widetilde \beta})f
=(J_{ \alpha}\otimes
J_{ \beta})
g$ reads
$$
f(x, y)=(J_{ \alpha}\otimes
J_{ \beta})g
=C_{\alpha}
C_{\beta}
\int_{\mathbb R^{2(n-1)}}
\frac 1{|x-z|^{2\alpha}
|y-w|^{2\beta}} g(z, w) dz dw.
$$
We perform the differentiation 
$\mathcal D_{\alpha, \beta, m}$
on this  identity 
and  find
$$
\mathcal D_{\alpha, \beta, m}f(x)=
C_{\alpha}
C_{\beta}
\int_{\br^{2(n-1)}}
S_{\alpha, \beta, m}(x; z, w) g(z, w) dz dw
=T_m g(x) =T_m
J_{\widetilde \alpha}\otimes
J_{\widetilde \beta} f(x),
$$
proving (\ref{eq:D-N-TJ}).

Finally the differentiation
under integral sign can be justified by  taking first
$\alpha, \beta << 0$
and $\alpha\notin \mathbb Z_-, \beta\notin \mathbb Z_-$,
with $\widetilde \alpha >>0, 
\widetilde \beta >>0$,  in which case  Lemma 2.1 implies that all integrals
involved are absolutely convergent. The meromorphic
continuation is proved in Lemma 3.2 and the
possible poles are read off from the definition of 
$\mathcal
D_{\alpha, \beta, m}$
and $\mathcal
M_{\alpha, \beta, m}$
 in  (\ref{eq:D-N}) respectively (\ref{eq:D-N-M}).
This finishes the proof.
\end{proof}

\section{Finitely many discrete components in the tensor product
$\pi_{\alpha}\otimes
\pi_{\beta}$, $G=SO_o(n, 1, \mathbb R)$}

We apply the intertwining operators $\mathcal D_m=\mathcal D_{\alpha,
  \beta, m}$
to
the study of appearance of discrete components
in the tensor product $\pi_{\alpha} \otimes
\pi_{\beta}    $ 
of complementary series. 
 For $\alpha, \beta\in (0, \rho)$ the tensor product 
$\pi_{\alpha} \otimes
\pi_{\beta}    $ in the non-compact picture
is the completion of $C^\infty_0(\mathbb R^{2(n-1)})$
with  norm
$$
\Vert 
f\Vert_{\alpha\otimes \beta}^2
:=\int_{\mathbb R^{2(n-1)}}
|\mathcal F f(\xi, \eta)
|^2 |\xi|^{n-2\alpha}
|\eta|^{n-2\beta} d\xi d\eta;
$$
cf.   (\ref{eq:Four-norm}).

\begin{theo+} Suppose $0<\alpha, \beta<\rho$ and $m\ge 0$ are integers.
If $\alpha +\beta +2m <\rho$ then the
intertwining operator $\mathcal D_{\alpha, \beta, m}$ is
a non-zero bounded intertwining operator
$\pi_{\alpha} \otimes
\pi_{\beta}
 \to  \pi_{\alpha+\beta +2m}$. Thus 
there are  $k$ discrete components $\pi_{\alpha+\beta +2m}$
 appears in the tensor product
$\pi_{\alpha} \otimes
\pi_{\beta}
$ for $m=0, 1, \cdots, k$ where $k$ is the maximal
integer such that $\alpha +\beta +2k <\rho$.
\end{theo+}

\begin{proof} Noticing that for $\alpha, \beta$ and $m$ as above we have
that the operator  $\mathcal D_m
=\mathcal D_{\alpha, \beta, m}$ is well-defined on smooth functions,
and 
$\pi_{\alpha}$, $\pi_{\beta }$ and $\pi_{\alpha+\beta +2m}$ are
unitary representations.
Recall also the notation $\tilde \alpha
=2\rho -\alpha=n-1-\alpha$
in \S2.3 and the unitary
norm (2.8). Let $f\in
C^\infty_0(\mathbb R^{2(n-1)})
\subset 
\pi_\alpha\otimes
\pi_\beta
$. 
We claim that 
$$ \Vert 
\mathcal D_m
f
\Vert_{\alpha+ \beta+2m}^2
\le C \Vert 
f\Vert_{\alpha\otimes \beta}^2.
$$
Thus $\mathcal D_m $
defines a non-zero intertwining operator from 
$\pi_{\alpha} \otimes
\pi_{\beta} $ into $
\pi_{\alpha+\beta +2m}$, proving our theorem.

Using Fourier inversion we have
$$
f(x, y)=(2\pi)^{-(n-1)}
\int_{\mathbb R^{2(n-1)}}
e^{i
\langle 
x, \xi
\rangle 
 +i
\langle 
 y, \eta
\rangle 
} 
\mathcal F f(\xi, \eta)
 d\xi d\eta.
$$
 We write
the differential operator $\mathcal E_{\alpha, \beta, m}$
in the proof of Theorem 3.1
as $Q(\mathcal L_x,  \mathcal L_y,
\nabla_x\cdot 
\nabla_y
)$
where $Q$ is a homogeneous polynomial of three variables
of degree $m$. 
Thus
$\mathcal D_mf(x)
=
Q(\mathcal L_x,  \mathcal L_y,
\nabla_x\cdot 
\nabla_y
)
f(x, y)|_{x=y}.
$ 
Its action on 
the inversion formula  results in
\begin{equation*}
  \begin{split}
\mathcal D_m
f(x)&=
(2\pi)^{-(n-1)}
\int_{\mathbb R^{2(n-1)}} 
e^{i
\langle  
x, \xi+ \eta
\rangle }
Q(-|\xi|^2,  -|\eta|^2,  -\langle
\xi, \eta \rangle)
\mathcal F f(\xi, \eta)
d\xi d\eta\\
&=(2\pi)^{-(n-1)}\int_{\mathbb R^{n-1}}
e^{i(x, \zeta)}
\int_{\mathbb R^{n-1}}
Q(-|\zeta-\eta|^2,  -|\eta|^2,  -\langle
\zeta-\eta, \eta \rangle)
\mathcal F f(\zeta-\eta, \eta)d\eta 
d\zeta.
  \end{split}
\end{equation*}
That is
\begin{equation*}
 \mathcal F(\mathcal D_m
f)(\zeta)=(2\pi)^{-(n-1)}
\int_{\mathbb R^{n-1}}
Q(-|\zeta-\eta|^2,  -|\eta|^2,  -\langle
\zeta-\eta, \eta \rangle)
\mathcal F f(\zeta-\eta, \eta)d\eta,
\end{equation*}
and furthermore 
\begin{equation*}
| \mathcal F(\mathcal D_m
f)(\zeta)|^2\le 
A(\zeta)
\int_{\mathbb R^{n-1}}
|\mathcal F f(\zeta-\eta, \eta)|^2
|\zeta-\eta|^{2\widetilde\alpha}
|\eta|^{2\widetilde\beta}
d\eta
\end{equation*}
with
$$
A(\zeta):
=(2\pi)^{-(n-1)}\int_{\mathbb R^{n-1}}
|Q(-|\zeta-\eta|^2,  -|\eta|^2,  -\langle
\zeta-\eta, \eta \rangle)|^2
|\zeta-\eta|^{-2\widetilde\alpha}
|\eta|^{-2\widetilde\beta} 
d\eta.
$$
To estimate the integral $A(\zeta)$ we write $\zeta=|\zeta|u$, $|u|=1$,
and perform a change of variables $\eta=|\zeta|v$. It is
$$
A(\zeta)
=(2\pi)^{-(n-1)}
|\zeta|^{4m-
2\widetilde\alpha
-2\widetilde\beta +(n-1)
}
\int_{\mathbb R^{n-1}}
|Q(-|u-v|^2,  -|v|^2,  -\langle
u-v, v \rangle)|^2
|u-v|^{-2\widetilde\alpha}
|v|^{-2\widetilde\beta} 
dv
$$
and the integral is convergent and uniformly bounded in $u, |u|=1$; indeed it is
locally integrable near  $v=0$, and $v=u$
for $2\widetilde\alpha, 2\widetilde\beta <n-1$ and
is integrable at infinity
for the integrand is dominated by
$$
C(1+|v|^2)^{ -(\widetilde\alpha
+\widetilde\beta-2m)}
$$
with $\widetilde\alpha
+\widetilde\beta-2m=n-1+ (n-1-\alpha-\beta-2m) <n-1$
for some constant $C$. 
Thus
$A(\xi) \le C|\zeta|^{4m-
2\widetilde\alpha
-2\widetilde\beta +(n-1)
}$, and 
$$
| \mathcal F(\mathcal D_m
f)(\zeta)|^2
|\zeta|^{-4m+
2\widetilde\alpha
+2\widetilde\beta -n
}\le C
\int_{\mathbb R^{n-1}}
|\mathcal F f(\zeta-\eta, \eta) |^2
|\zeta-\eta|^{2\widetilde\alpha}
|\eta|^{2\widetilde\beta}
d\eta,
$$
and its integration  over $\zeta$ gives
\begin{equation*}
  \begin{split}
\int_{\mathbb R^{n-1}}
| \mathcal F(\mathcal D_m
f)(\zeta)|^2
|\zeta|^{-4m+
2\widetilde\alpha
+2\widetilde\beta -n
}d\zeta &
\le C \int_{\mathbb R^{n-1}}
\int_{\mathbb R^{n-1}}
|\mathcal F f(\zeta-\eta, \eta) |^2
|\zeta-\eta|^{2\widetilde\alpha}
|\eta|^{2\widetilde\beta}
d\eta\\
&= C\Vert 
f\Vert_{\alpha\otimes \beta}^2    
  \end{split}
\end{equation*}
whereas the LHS is precisely
$ \Vert 
\mathcal D_mf\Vert_{\alpha+ \beta+2m}^2$.
This finishes the proof.
\end{proof}

When $n=2$ then   $m=0$ and the  theorem
states
that 
$\pi_{\alpha+\beta}$ appears in the tensor product
$\pi_{\alpha} \otimes
\pi_{\beta}$ if $\alpha+\beta <1$. This has been
 proved earlier in \cite{Repka78}.

\section{The appearance of 
one component $\pi_{\alpha+\beta}$ in $\pi_{\alpha}\otimes
\pi_{\beta}$ for  $G=SU(n, 1),
Sp(n, 1)$}

We treat now the other rank one classical groups. The following
result follows straightforward
from Theorem 2.1.

\begin{theo+}  Let $G=SU(n, 1)$ and $Sp(n, 1)$,
$\pi_{\alpha}$ and $\pi_{\beta}$
be the complementary series for $\alpha, \beta$ as in
(\ref{mu-fixed}),
$0<\alpha, \beta <\rho=n$ and respectively
$2<\alpha, \beta <\rho=2n-1$.
Then the complementary series
$(\pi_{\alpha+\beta}, G)$ of
$G$
appears discretely in the tensor product
$\pi_{\alpha}\otimes \pi_{\beta}$ 
if
 $$
\alpha +\beta <
 \begin{cases} 
n,  & \mathbb F=\bc, \\
2n-1,\quad & \mathbb F=\bh.
 \end{cases}
$$
\end{theo+}

\begin{proof} We prove the case for $G=SU(n, 1)$ and the same methods applies also
to $G=Sp(n, 1)$.  We consider
the diagonal imbedding of $G$
in $G_1=SU(n, 1)
\times SU(n, 1)
$. It follows from Theorem 2.1
that for $\alpha, \beta\in (0, \rho)$
the complementary series $\pi_\alpha$
and $\pi_\beta$ appear in $\tau_{\frac \alpha
2 }$
and $\tau_{\frac \beta 2
}$, respectively.
Now $\tau_{\nu}$ of 
$G_1
$
is the tensor product
$\lambda_\nu \otimes \overline
{\lambda_{\nu}}$ on 
$\mathcal H_\nu
\otimes \overline
{\mathcal H_\nu}
$
where $\mathcal H_\nu$ is the space
of holomorphic functions on the unit ball
$B^n$ with the reproducing kernel $(1-\langle z, w\rangle)^{-\nu}$.
If $\alpha +\beta<n$
then 
$\pi_\alpha$
 appears in $\tau_{\frac \alpha 2}$, so does
$\pi_\beta$ in
$\tau_{\frac \beta 2}$.
 The tensor product
 $\tau_{\frac \alpha 2
}\otimes
\tau_{\frac \beta 2}$ is now
$$
H:=(\mathcal H_{\frac \alpha 2}
\otimes \overline
{\mathcal H_{\frac \alpha 2}} )
\otimes
(\mathcal H_{\frac \beta 2}
\otimes \overline
{\mathcal H_{\frac \beta 2}} ).
$$
Its restriction to $G$
is
$$H=(\mathcal H_{\frac \alpha 2}
\otimes
\mathcal H_{\frac \beta 2})\otimes
 \overline
{(\mathcal H_{\frac \alpha 2}
\otimes\mathcal H_{\frac \beta 2})
}.
$$
However the tensor product 
$\mathcal H_{\frac \alpha 2}
\otimes
\mathcal H_{\frac \beta 2}$ 
of two holomorphic representation is
decomposed discretely under $G$ and contains
a component 
$\mathcal H_{\frac{\alpha+\beta }2}
$; see e.~g. \cite{Repka79, p-z-tensor}  and references therein. Thus $H$ contains a discrete component 
$
\mathcal  H_{\frac{\alpha+\beta }2}
\otimes
 \overline
{\mathcal 
 H_{\frac{\alpha+\beta }2}
}.
$
We use again Theorem 2.1 and deduce that
this space has a discrete component $(\pi_{\alpha +\beta}, G)$.
\end{proof}

\def\cprime{$'$} \newcommand{\noopsort}[1]{} \newcommand{\printfirst}[2]{#1}
  \newcommand{\singleletter}[1]{#1} \newcommand{\switchargs}[2]{#2#1}
  \def\cprime{$'$} \def\cprime{$'$} \def\cprime{$'$}
\providecommand{\bysame}{\leavevmode\hbox to3em{\hrulefill}\thinspace}
\providecommand{\MR}{\relax\ifhmode\unskip\space\fi MR }
\providecommand{\MRhref}[2]{%
  \href{http://www.ams.org/mathscinet-getitem?mr=#1}{#2}
}
\providecommand{\href}[2]{#2}


\begin{thebibliography}{10}

\bibitem{Asmuth-Repka}
C.~Asmuth and J.~Repka, \emph{Tensor products for { ${\rm SL}_{2}(\mathcal K)$
  }. {I}. {C}omplementary series and the special representation}, Pacific J.
  Math. \textbf{97} (1981), no.~2, 271--282. \MR{641157 (83c:22023)}

\bibitem{BKZ} S. Ben Sa{\" \i}d, K. Koufany and G. Zhang,
\emph
{Invariant trilinear forms on  spherical principal series of  real-rank
one semisimple  Lie groups,} Internat. J. Math. 25 (2014), no. 3, 
1450017 (35 pages).


\bibitem{Burger-Li-Sarnak}
M.~Burger, J.-S. Li, and P.~Sarnak, \emph{Ramanujan duals and automorphic
  spectrum}, Bull. Amer. Math. Soc. (N.S.) \textbf{26} (1992), no.~2, 253--257.
  \MR{1118700 (92h:22023)}

\bibitem{Clerc-2011}
J.-L. Clerc, \emph{Singular conformally invariant trilinear forms and covariant
  differential operators on the sphere}, arXiv:1102.1861.

\bibitem{Clerc-Beckmann}
J.-L. Clerc and R.~Beckmann, \emph{Singular conformally invariant trilinear
  forms and generalized {Rankin-Cohen} operators}, arXiv:1104.3461.

\bibitem{Clerc-Orsted}
J.-L. Clerc and B.~\O{}rsted, \emph{Conformally invariant trilinear forms on
  the sphere}, arXiv:1001.2851.

\bibitem{Clozel-ias-park-lec}
L.~Clozel, \emph{Spectral theory of automorphic forms}, Automorphic forms and
  applications, IAS/Park City Math. Ser., vol.~12, Amer. Math. Soc.,
  Providence, RI, 2007, pp.~43--93. \MR{2331344 (2008i:11069)}

\bibitem{Connes-Mosc-2}
A.~Connes and H.~Moscovici, \emph{Rankin-{C}ohen brackets and the {H}opf
  algebra of transverse geometry}, Mosc. Math. J. \textbf{4} (2004), no.~1,
  111--130.

\bibitem{CDKR-jga}
M.~Cowling, A.~Dooley, A.~Kor{\'a}nyi, and F.~Ricci, \emph{An
  approach to symmetric spaces of rank one via groups of {H}eisenberg type}, J.
  Geom. Anal. \textbf{8} (1998), no.~2, 199--237.

\bibitem{He2}
S.~Helgason, \emph{Groups and geometric analysis}, Academic Press, New York,
  London, 1984.

\bibitem
{CKOP}
B.~\O{}rsted J.-L.~Clerc, T.~Kobayashi and M.~Pevzner, 
\emph{
Generalized
  {Bernstein--Reznikov} integrals}, Math. Ann. \textbf{349} (2011), 395--431.

\bibitem{Dijk-Hille-jfa}
G. van Dijk and S. C. Hille,
\emph{
Canonical representations
related to hyperbolic spaces,
}
J. Funct. Anal.,
\textbf{147}  (1997), 109-139.

\bibitem{FK}
J.~Faraut and A.~Koranyi, \emph{Function spaces and reproducing kernels on
  bounded symmetric domains}, J. Funct. Anal. \textbf{88} (1990), 64--89.


\bibitem{Johnson-Wallach}
K.~D.~Johnson and N.~R.~Wallach,
 \emph{
Composition series and
  intertwining operators for the spherical principal series. {I}}, Trans. Amer.
  Math. Soc. \textbf{229} (1977), 137--173. \MR{MR0447483 (56 \#5794)}


\bibitem{Kn-book}
A.~Knapp, \emph{Representation theory of semisimple groups}, Princeton
  University Press, Princeton, New Jersey, 1986.



\bibitem{Kob-Pev}
 T. Kobayashi and  M. Pevzner,
 \emph{
Rankin-Cohen Operators for Symmetric Pairs,}
preprint, arXiv:1301.2111.

\bibitem{Molchanov}
V.~F.~Molchanov,
 \emph{
Canonical representations and overgroups,}
Lie groups and symmetric spaces, 213-224, 
Amer. Math. Soc. Transl. Ser. 2, 210, Adv. Math. Sci., 54, Amer. Math. Soc., Providence, RI, 2003.


\bibitem{Mollers-Oshima}
J.~M\"o{}llers and Y.~Oshima,
 \emph{
Restriction of most degenerate representations of O(1,N) with respect to symmetric pairs,}
J. Math. Sci. Univ. Tokyo \textbf{22} (2015), no. 1, 279-338. 


\bibitem{Naimark-1961}
M.~A. Na{i}mark, \emph{Decomposition of a tensor product of irreducible
  representations of the proper {L}orentz group into irreducible
  representations. {III}}, Trudy Moskov. Mat. Ob\v s\v c. \textbf{10} (1961),
  181--216. \MR{0148800 (26 \#6304b)}

\bibitem{Neretin-1995}
Yu.~A. Neretin and G.~I. Ol{\cprime}shanski{i}, \emph{Boundary values of
  holomorphic functions, singular unitary representations of the groups {${\rm
  O}(p,q)$} and their limits as {$q\to\infty$}}, 
J. Math. Sci. (New York) \textbf{87} (1997), no. 6, 3983 4035 



\bibitem{Ovsienko-Redou}
V.~Ovsienko and P.~Redou, \emph{Generalized transvectants-{R}ankin-{C}ohen
  brackets}, Lett. Math. Phys. 
\textbf{63} (2003), no.~1, 19--28. \MR{1967533
  (2004c:58077)}

\bibitem{p-z-tensor}
L.~Peng and G.~Zhang, \emph{Tensor products of holomorphic
  representations and bilinear differential operators}, J. Funct. Anal.
  \textbf{210} (2004), no.~1, 171--192.

\bibitem{Pukanszky62}
L.~Puk{\'a}nszky, \emph{On the {K}ronecker products of irreducible
  representations of the {$2\times 2$} real unimodular group. {I}}, Trans.
  Amer. Math. Soc. \textbf{100} (1961), 116--152. \MR{0172962 (30 \#3177)}

\bibitem{Repka78}
J.~Repka, \emph{Tensor products of unitary representations of ${SL_2(\bf R)}$},
  Amer. J Math. \textbf{100} (1978), 930--932.

\bibitem{Repka79}
J.~Repka, \emph{Tensor products of holomorphic discrete series
  representations}, Can. J. Math. \textbf{31} (1979), 836--844.

\bibitem
{Speh-Venk-2}
B.~Speh and T.~N. Venkataramana, 
\emph{
Discrete components of some
  complementary series representations}, Indian J. Pure Appl. Math. \textbf{41}
  (2010), no.~1, 145--151. \MR{2650105 (2011h:22014)}

\bibitem{Speh-Zhang}
B. Speh and  G. Zhang,
\emph{ Restriction to symmetric subgroups of unitary representations
  of rank one semisimple Lie groups}, 
Math. Zeit.
\textbf{283}, 629-647.

\bibitem{Vogan-Wallach}
D.~A. Vogan, Jr. and N.~R. Wallach, 
\emph{
Intertwining operators for real
  reductive groups}, Adv. Math. \textbf{82} (1990), no.~2, 203--243.
  \MR{1063958 (91h:22022)}

\bibitem{Zhang-jfa}
 G. Zhang,
\emph{
 Discrete components in restriction of unitary representations of rank one semisimple Lie groups},
 J. Funct. Anal. \textbf{269} (2015), no. 12, 3689-3713.


\end{thebibliography}

\end{document}